\documentstyle[11pt,amstex,amsfonts,amssymb]{amsart}
\newcommand{\rest}{{\,|\grave{}\,}}

\newcommand{\V}{{\bf V}} 
\newcommand{\conc}{{}^\frown\!}
\newcommand{\lh}{\ell g\/}

\newcommand{\sk}{{\bf S}_\kappa}
\newcommand{\bkzo}{{\square^\kappa[0,1]}}
\newcommand{\QED}{\hfill\vrule width 6pt height 6pt depth 0pt 
\vspace{0.1in}} 
\newcommand{\Proof}{\noindent{\sc Proof} \hspace{0.2in}} 
\newcommand{\Hchi}{{\cal H}(\chi)}

\newcommand{\ZFC}{{\rm ZFC}}
\newcommand{\MA}{{\rm MA}}
\newcommand{\con}{{2^{\aleph_0}}}
\newcommand{\can}{{}^{\textstyle \omega}2}
\newcommand{\fso}{{}^{\textstyle \omega>}\omega}
\newcommand{\cl}{{\rm cl}}
\newcommand{\gr}{{\frak r}}
\newtheorem{theorem}{Theorem}
\newtheorem{fact}[theorem]{Fact}
\newtheorem{lemma}[theorem]{Lemma} 
\newtheorem{proposition}[theorem]{Proposition} 
\newtheorem{corollary}[theorem]{Corollary} 

\theoremstyle{definition}

\theoremstyle{remark}

\title{A result related to the problem CN of Fremlin}  

\author[O. Kolman]{\bf \uppercase {O. Kolman}}

\address{ 2 West Eaton Place, London SW1X 8LS, U.K.}

\author[S. Shelah]{\uppercase {\bf S. Shelah}}
\thanks{The research partially supported by
``Basic Research Foundation'' of the Israel Academy of Sciences and
Humanities. Publication 661}

\address{Institute of Mathematics\\
The Hebrew University\\
Jerusalem 91904, Israel\\
and Rutgers University\\
Mathematics Department\\
New Brunswick, NJ 08854, USA
}
\email{shelah@@math.huji.ac.il}

\date{\today}

\begin{document}
\maketitle

\bigskip
\bigskip
\bigskip
{\it Abstract.}
We show that the set of injective functions from any uncountable cardinal less
than the continuum into the real numbers is of second category in the box
product topology.
\bigskip
\bigskip
\bigskip

In this paper, we use a definability argument to resolve under mild
set--theoretic assumptions a problem about injective functions in the box
product topology. Suppose $\kappa$ is a cardinal; let $\sk$ be the set of
injective functions from $\kappa$ into the closed unit interval $[0,1]$, and
let $\bkzo$ be the product of $\kappa$ copies of $[0,1]$ equipped with the box
product topology (basic open sets are just products of basic open subsets of
$[0,1]$: is $\sk$ of first category in $\bkzo$. We shall prove:

\begin{theorem}
\label{main}
Suppose that $\aleph_1\leq\kappa<\con$. Then $\sk$ is of second category in
$\bkzo$. 
\end{theorem}

David Fremlin asks \cite[Problem CN]{Fe94} whether ${\bf S}_{\aleph_1}$ is
co-meagre and whether ${\bf S}_{\con}$ is of the second category. We do not
know the answer to these questions, but feel that Theorem \ref{main}
represents some progress on the former. As regards the latter question, we
note that it is easy to prove that if $\MA$ holds, then every countable
intersection of dense open sets in $\square^\con[0,1]$ contains functions
which are injective on a closed unbounded subset of $\con$ (see Corollary
\ref{corfinal}).  

In proving Theorem \ref{main}, we make use (often tacitly) of some standard
results about elementary submodels of $\Hchi$, the set of all sets
hereditarily of cardinality less than $\chi$. For reader's convenience we
record these next.

\begin{fact}
\begin{enumerate}
\item The cardinals $\omega$ and $\omega_1$ belong to $\Hchi$ for
$\chi>\omega_1$.
\item If $\chi>\omega_1$ is regular, then $\Hchi$ is a model of all the axioms
of $ZFC$ except possibly the power set axiom.
\item If $\chi>\omega_1$ is singular, then $\Hchi$ is a model of all the
axioms of $\ZFC$ except possibly the power set, union and replacement axioms. 
\QED 
\end{enumerate}
\end{fact}

\begin{lemma}
\label{lemma}
Suppose that $\chi>\omega_1$ is regular and $N$ is an elementary submodel of
$\Hchi$. Then:
\begin{enumerate}
\item If $a\in\Hchi$ is definable with parameters from $N$ (i.e., there is a
formula $\varphi(x)$ with one free variable $x$ and possibly parameters in $N$
such that $a$ is the unique element which satisfies $\varphi(x)$ in
$(\Hchi,\in)$), then $a\in N$. 
\item The ordinals $\omega$ and $\omega_1$ belong to $N$ and $\omega\subseteq
N$. 
\item If $a,A,B,f\in N$ and (in $\V$) $f$ is a function from $A$ to $B$ then
$f(a)\in N$.
\item For every $\alpha\in \omega_1\cup\{\omega_1\}$, $\alpha\cap N$ is an
ordinal.
\item If $\kappa\in N$ and $\{A_\alpha:\alpha<\kappa\}\in N$, then $(\forall
\alpha\in\kappa\cap N)(A_\alpha\in N)$. \QED
\end{enumerate}
\end{lemma}

The proofs of these well-known facts can be found in many places, for example,
in \cite{EM} or in the appendix of \cite{HS}. 

Since the proof of Theorem \ref{main} involves some notation, we describe in
intuitive terms how it works. First we shall fix a family of ``rational
boxes'' which are products of open subintervals of $[0,1]$ having rational
endpoints. For a given countable family of dense open subsets $G_n$ of
$\bkzo$, for each $n$, we define a function $f_n$ on rational boxes such that
$f_n(B)\subseteq G_n$ and the diameters of open sets used to define the box
$f_n(B)$ decrease as $n$ increases. For each real $\tau$, we define using the
$f_n$'s an element $x_\tau\in\bigcap\limits_{n<\omega} G_n$. We then take an
elementary submodel $N$ of cardinality $\kappa$ (of $\Hchi$ for a large enough
$\chi$) containing the $f_n$'s, $G_n$'s, all the ordinals up to $\kappa+1$,
and whatever else is necessary. Since $\kappa<\con$, there is a real
$\tau\in\can\setminus N$. We complete the proof by showing that $x_\tau\in
\sk$, and this is done by demonstrating that if $x_\tau$ is not injective,
then $x_\tau$ is definable in $N$ and hence belongs to $N$ -- a
contradiction. So the essence of the argument lies in connecting
non-injectivity and definability. 
\medskip

\noindent{\sc Proof of Theorem \ref{main}} \hspace{0.2in} Let $\{G_n:n\in
\omega\}$ be a countable family of dense open sets in $\bkzo$. We must show
that $\sk\cap\bigcap\limits_{n<\omega} G_n\neq\emptyset$.

First we set up some notation.

Fix a canonical family of non-empty open intervals $I_\rho$ for $\rho\in\fso$
as follows. Let $I_{\langle\rangle}=(0,1)$. Suppose that $I_\rho$ has been
defined, $I_\rho$ a non-empty open subinterval of $[0,1]$, and the length
$\lh(\rho)$ of the sequence $\rho$ is $n$. Choose $2^{2^n}$ disjoint open
subintervals of $I_\rho$ of equal length, say $I_{\rho\conc\langle k\rangle}$,
$0\leq k<2^{2^n}$ such that $I_\rho\setminus \bigcup\limits_{k<2^{2^n}}
I_{\rho\conc\langle k\rangle}$ is finite. This completes the definition of the
family $\{I_\rho:\rho\in\fso\}$. 

Next for each $n\in\omega$, we choose a function $f_n$ defined on the family
of non-empty open boxes $B=\prod\limits_{\iota<\kappa}B_\iota$, where each
$B_\iota$ is a non-empty open subinterval of $[0,1]$, as follows:
$f_n(B)=\prod\limits_{\iota<\kappa} A_{\iota,n}$, where for all
$\iota<\kappa$:
\begin{enumerate}
\item $A_{\iota,n}\in\{I_\rho:\lh(\rho)>n\}$,
\item $\cl(A_{\iota,n})\subseteq B_n$,
\item $f_n(B)\subseteq G_n$.
\end{enumerate}
There is no problem in choosing $f_n$ as above, since $G_n$ is dense open and
for each $\iota<\kappa$ we can take $\lh(\rho)$ large enough to ensure that
$\cl(I_\rho)\subseteq B_\iota$.

Finally, fix $\kappa$ reals $\{\eta_i:i<\kappa\}$, $\eta_i\in\can$, $\eta_i
\neq\eta_j$ for $i<j<\kappa$.

We associate with each real $\tau$ an element $x_\tau\in\bigcap\limits_{n<
\omega} G_n$. Define by induction on $n$, for every $i<\kappa$, a non-empty
open subinterval $C_{i,\tau\rest n}\subseteq [0,1]$ and an open box $C_{\tau
\rest n}=\prod\limits_{i<\kappa} C_{i,\tau\rest n}$ as follows. Let $C_{i,\tau
\rest 0}=(0,1)$ (for $i<\kappa$). Suppose that $C_{i,\tau\rest n}$ has been
defined (for $i<\kappa$) and is a non-empty open subinterval of $[0,1]$. By
(1), $f_n(C_{\tau\rest n})=\prod\limits_{i<\kappa} I_{\rho_i}$, for some
$\rho_i\in\fso$ such that $\lh(\rho_i)>n$. Let 
\[k_i=2\cdot |\{\eta\in {}^{\textstyle n}2: \eta\leq_{\ell ex}\eta_i\rest n\}|
+\tau(n).\]
Note that (trivially) $k_i<2^{2^{\lh(\rho_1)}}$, so $I_{\rho_i\conc\langle k_i
\rangle}$ is a non-empty open subinterval. Put $C_{i,\tau\rest n+1}= I_{\rho_i
\conc\langle k_i\rangle}$. By (1), (2) and (3), there exists a unique element
\[x_\tau\in\bigcap_{n<\omega} C_{\tau\rest n}\subseteq \bigcap_{n<\omega}
G_n.\]
So to complete the proof, it will suffice to choose a real $\tau$ such that
$x_\tau\in\sk$. 

Let $N\prec\Hchi$ be an elementary submodel for $\chi$ regular large enough
($(\beth_\omega)^+$ will do) such that
\begin{enumerate}
\item $\{f_n:n\in\omega\}\subseteq N$, $\kappa+1\subseteq N$, $\{\eta_i:\i<
\kappa\}\subseteq N$, $\{G_n:n\in\omega\}\subseteq N$, $\{I_\rho: \rho\in\fso
\}\subseteq N$,
\item $|N|=\kappa$.
\end{enumerate}
Since $|N|=\kappa<\con$, there exists a real $\tau\in\can\setminus N$. We
complete the proof by proving that $x_\tau$ is an injective function from
$\kappa$ into $[0,1]$.

Suppose that $x_\tau$ is not injective: so there are $i<j<\kappa$ such that
$x_\tau(i)=x_\tau(j)$. We derive a contradiction by showing that $\tau$ is
definable using parameters in $N$ and hence $\tau\in N$ by Lemma \ref{lemma}. 

Let $m_0=\min\{n:\eta_i(n)\neq \eta_j(n)\}$ (we can calculate $m_0$ in $N$
since $\{\eta_i:i<\kappa\}\subseteq N$ and $i,j\in N$). It suffices to show
that we can define $\tau\rest n$ in $N$ for every $n>m_0$. We prove this by
induction on $n$. Suppose that we have defined $\tau\rest n$, $n>m_0$. We show
how to calculate $\tau(n)$ in $N$, thereby defining $\tau\rest (n+1)$. Note
that $C_{i,\sigma}$ and $C_\sigma$ are definable in $N$ (for $i<\kappa$ and
$\sigma\in {}^{\textstyle \omega>}2$) and hence $f_n(C_{\tau\rest n})=
\prod\limits_{\alpha<\kappa} I_{\rho_\alpha}\in N$. Consider $I_{\rho_i}$ and
$I_{\rho_j}$: in $N$, $I_{\rho_i}\cap I_{\rho_j}\neq\emptyset$ since
$x_\tau(i)=x_\tau(i)$ and $N\prec\Hchi$. Also in $N$, $\rho_i\leq_{\ell ex}
\rho_j$ or $\rho_j\leq_{\ell ex}\rho_i$, so without loss of generality,
$\rho_i\leq_{\ell ex}\rho_j$. In fact, $\rho_i<_{\ell ex}\rho_j$, for if
$\rho_i=\rho_j$, then $\eta_i\rest n=\eta_j\rest n$, contradicting $n>m_0$.
Since $I_{\rho_i}\cap I_{\rho_j}\neq\emptyset$, it follows that $\lh(\rho_i)<
\lh(\rho_j)$ and $\rho_j\rest\lh(\rho_i)=\rho_i$, and so there is a unique
natural number $k^*$ such that $\rho_j\rest (\lh(\rho_i)+1)=\rho_i\conc\langle
k^*\rangle$. As we can compute $k^*$ and $2\cdot |\{\eta\in{}^{\textstyle n}2:
\eta\leq_{\ell ex}\eta_i\rest n\}|$ inside $N$, it will suffice to show (in
$\Hchi$) that: 
\begin{enumerate}
\item[$(*)$]\quad $\tau(n)=k^*-2\cdot |\{\eta\in {}^{\textstyle n}2: \eta
\leq_{\ell ex}\eta_i\rest n\}|$.
\end{enumerate}
To see this, notice that $C_{i,\tau\rest (n+1)}=I_{\rho_i\conc\langle k^*
\rangle}$. Why? If not, then $C_{i,\tau\rest (n+1)}= I_{\rho_i\conc\langle
\ell \rangle}$ for some $\ell\neq k^*$, and so $I_{\rho_i\conc \langle k^*
\rangle}\cap I_{\rho_i\conc\langle\ell\rangle}=\emptyset$, and hence 
\begin{enumerate}
\item[$(**)$]\quad $I_{\rho_j}\cap I_{\rho_i\conc\langle\ell\rangle}=
\emptyset$ 
\end{enumerate}
(as $I_{\rho_j}\subseteq I_{\rho_i\conc\langle k^*\rangle}$). However,
$x_\tau(i)\in C_{i,\tau\rest (n+1)}=I_{\rho_i\conc\langle\ell\rangle}$, and
$x_\tau(i)=x_\tau(j)\in I_{\rho_j}$, so $I_{\rho_j}\cap I_{\rho_i\conc\langle
\ell\rangle}\neq\emptyset$, contradicting $(**)$. By the definition of
$C_{i,\tau\rest(n+1)}=I_{\rho_i\conc\langle k_i\rangle}$, it is now immediate
that $k_i=k^*$. Recalling that  
\[k_i=2\cdot |\{\eta\in {}^{\textstyle n}2:\eta\leq_{\ell ex}\eta_i\rest n\}|
+\tau(n),\]
we obtain $(*)$, and so $\tau\rest (n+1)$ is definable in $N$.

Thus, $\tau=\bigcup\limits_{n<\omega}(\tau\rest n)$ is definable using
parameters in $N$ and hence belongs to $N$ -- a contradiction. It follows that
$x_\tau$ is injective. This completes the proof of the theorem. \QED

\begin{proposition}
Suppose that $\con$ is regular and there is an enumeration $\{\gr_\alpha:
\alpha<\con\}$ of the real numbers such that for every $\beta$, $\{\gr_\alpha:
\alpha<\beta\}$ is of the first category. Then every countable intersection of
dense open sets in $\square^\con[0,1]$ contains functions which are injective
on a closed unbounded set of $\con$.
\end{proposition}

\Proof Let $X_\beta=\{\gr_\alpha:\alpha<\beta\}\subseteq\bigcup\limits_{n<
\omega}F_{\beta,n}$, where $F_{\beta,n}$ is closed nowhere dense in $[0,1]$,
and let $G_{\beta,n}=\can\setminus F_{\beta,n}$. So $G_n=\prod\limits_{\alpha<
\con}G_{\alpha,n}$ is dense and open. It is enough to show that every $f\in
\bigcap\limits_{n<\omega} G_n$ is injective on some club. Suppose that $f\in
\bigcap\limits_{n<\omega} G_n$. For each $\alpha<\con$, let 
\[\beta_\alpha\stackrel{\rm def}{=}\min\{\beta:f(\alpha)=\gr_\beta\}.\]
The interval of ordinals $(\beta_\alpha,\con)$ is a club, and hence the
diagonal intersection  
\[C=\mathop{\triangle}_\alpha(\beta_\alpha,\con)=\{\beta<\con: (\forall\alpha
<\beta)(\beta\in (\beta_\alpha,\con))\}\]
is a club. It is trivial to verify that if $\beta\in C$ and $\alpha<\beta$
then $f(\alpha)\neq f(\beta)$. \QED

\begin{corollary}
\label{corfinal}
Martin's Axiom implies that every countable intersection of dense open subsets
in $\square^\con[0,1]$ contains functions which are injective on some
club. \QED 
\end{corollary}

\end{document}